\documentclass{amsart}
\usepackage{graphicx}
\usepackage{amssymb}

\newtheorem{thm}{Theorem}[section]
\newtheorem{cor}[thm]{Corollary}

\theoremstyle{definition}

\theoremstyle{remark}
\newtheorem{rem}[thm]{Remark}

\begin{document}

\title[Singularities of complex-valued solutions]{Singularities of complex-valued solutions to linear parabolic equations}
\author{Connor Mooney}
\address{Department of Mathematics, ETH Z\"{u}rich}
\email{\tt connor.mooney@math.ethz.ch}


\begin{abstract}
We construct examples of blowup from smooth data for complex-valued solutions to linear uniformly parabolic equations in dimension $n \geq 2$, which are exactly as irregular as parabolic energy estimates allow.
\end{abstract}
\maketitle
\section{Introduction}
In this paper we consider linear uniformly parabolic equations of the form
\begin{equation}\label{Equation}
 u_t - \text{div}(A(x,\,t)\nabla u) = 0.
\end{equation}
Here $u : \mathbb{R}^{n+1} \rightarrow \mathbb{C}$, and the coefficients are bounded measurable, complex-valued functions satisfying
\begin{equation}\label{Ellipticity}
 \text{Re}(A_{kl}(x,\,t)p_k\overline{p}_l) \geq \lambda |p|^2, \quad |A(x,\,t)(p)|^2 \leq \Lambda^2 |p|^2
\end{equation}
for some constants $\lambda,\, \Lambda > 0$, and for all $(x,\,t) \in \mathbb{R}^{n+1}$ and $p \in \mathbb{C}^n$. By a solution we mean that $u \in L^2_{loc,\,t}(H^1_{loc,\,x})$ solves (\ref{Equation})
in the sense of distributions. We note that (\ref{Equation}) can be viewed as a uniformly parabolic system of the form
\begin{equation}\label{System}
 \partial_tv^{\alpha} - \partial_k(B^{kl}_{\alpha\beta}(x,\,t)v^{\beta}_l) = 0, \quad 1 \leq k,\,l \leq n, \quad 1 \leq \alpha,\,\beta \leq 2.
\end{equation}
Here $u = v^1 + iv^2$ and $B_{11} = B_{22} = \text{Re}(A),\, B_{12} = -B_{21} = -\text{Im}(A)$.

We briefly discuss the elliptic case
\begin{equation}\label{EllipticEquation}
 \text{div}(A(x)\nabla u) = 0
\end{equation}
in $\mathbb{R}^n$. Solutions to (\ref{EllipticEquation}) are $C^{\alpha}$ when $n = 2$ by work of Morrey \cite{Mo}. Real-valued solutions are $C^{\alpha}$
by fundamental work of De Giorgi \cite{DG1} and Nash \cite{Na}. 
There are classical counterexamples to continuity for solutions to elliptic systems when $n \geq 3$ (see \cite{DG2}, \cite{GM}, \cite{Ma}). Discontinuous solutions to (\ref{EllipticEquation})
were first constructed in dimension $n \geq 5$ \cite{MNP}, and later in dimension $n \geq 3$ \cite{F}. In general, the best regularity we have for (\ref{EllipticEquation})
is $u \in W^{1,\,2 + \delta}_{loc}$ for some $\delta(n,\,\lambda,\,\Lambda) > 0$ (see \cite{Gi}), which is only slightly better than the energy class of the solutions. 
In fact, for each $\gamma > 2$ there are solutions to (\ref{EllipticEquation}) that are not in $W^{1,\,\gamma}_{loc}$ (see \cite{F}).

Interestingly, the parabolic problem (\ref{Equation}) has resisted a similar understanding. Real-valued solutions are $C^{\alpha}$ \cite{Na}. In general
we have the higher-integrability results $\nabla u \in L^{2 + \delta}_{loc}$ and $u \in L^{\infty}_{loc,\, t}(L^{2+\delta}_{loc, \, x})$ for some $\delta(n,\,\lambda,\,\Lambda) > 0$ (see \cite{St}, \cite{NS}).
There are also examples of discontinuity from smooth data when $n \geq 3$ (\cite{FM}, and \cite{SJM}, \cite{SJ2} for more general systems). 
However, the examples are in $L^{\infty}_{loc,\,t}(W^{1,\,2+\delta}_{loc,\,x})$, and are thus significantly more regular than the higher-integrability results predict.
When $n = 2$ the known results don't imply continuity of solutions (unlike the elliptic case),
which remained open for some time (see e.g. \cite{SJM}, \cite{JS}, \cite{SJ1}, \cite{SJ2}). We recently settled this problem with a counterexample \cite{M1}. 
Still, the example in \cite{M1} is barely irregular enough to develop a discontinuity (it is e.g. in $L^{\infty}_{loc,\,t}(L^{p}_{loc,\,x})$ for $p$ large), so the regularity gap between theory and examples remained large.

The purpose of this paper is to complete the picture for (\ref{Equation}) by constructing solutions
in dimension $n \geq 2$ that are exactly as irregular as the parabolic higher-integrability results allow. We also prove some Liouville theorems which explain why previous approaches only produced ``elliptic" discontinuities.
Our results connect the regularity problem for (\ref{Equation}) in $\mathbb{R}^{n+1},$ in parabolic geometry, to that for the elliptic equation (\ref{EllipticEquation}) in $\mathbb{R}^{n+2}$.
We make this connection precise in the next section.

\section{Results}
In this section we state our results. We will deal with ``spiraling'' self-similar solutions to (\ref{Equation}) of the form
\begin{equation}\label{Ansatz}
u(x,\,t) = (-t)^{-\frac{\mu}{2}}\,e^{-\frac{i}{2}\log(-t)}\, w\left(\frac{x}{(-t)^{1/2}}\right).
\end{equation}
These are invariant under $u \rightarrow \lambda^{\mu}e^{i\log \lambda} u(\lambda x,\, \lambda^2t).$ We obtain a solution
to (\ref{Equation}) on $\mathbb{R}^n \times (-\infty,\,0)$ with coefficients $A(x/(-t)^{1/2})$ if $w$ solves the elliptic equation
\begin{equation}\label{SelfSimEquation}
 \text{div}(A(x)\nabla w) = \frac{1}{2}(iw + \mu w + x \cdot \nabla w)
\end{equation}
on $\mathbb{R}^n$, and $A$ satisfies (\ref{Ellipticity}) for some $\lambda,\,\Lambda > 0$. 
Furthermore, the solution defined by (\ref{Ansatz}) is smooth up to $t = 0$ away from $x = 0$ and develops a ``spiraling $-\mu$-homogeneous'' discontinuity at $t = 0$ provided $\mu \geq 0$ and
\begin{equation}\label{Asymptotics}
 w = |x|^{-\mu}g(x/|x|)e^{-i\log |x|}(1 + \mathcal{E}(|x|^{-2})) \,\, \text{ on } \,\, \mathbb{R}^n \backslash B_1.
\end{equation}
Here $g \in C^{\infty}(S^{n-1})$ and $\mathcal{E}$ is a smooth function with $\mathcal{E}(0) = 0$. We can extend the solution to positive times e.g. by solving the heat equation with initial data $u(x,\,0) := |x|^{-\mu}g(x/|x|)e^{-i\log |x|}$,
provided $\mu < n$.

Our first result is:
\begin{thm}\label{Counterexample}
 If $n \geq 2$ and $0 \leq 2\mu < n$, then there exists a nontrivial solution to (\ref{SelfSimEquation}) on $\mathbb{R}^n$ that satisfies (\ref{Asymptotics}).
\end{thm}
\noindent By taking $\mu$ arbitrarily close to $\frac{n}{2}$ we obtain as a consequence:
\begin{cor}\label{Optimality}
 For all $n \geq 2$ and $\delta > 0$, there exists a solution to (\ref{Equation}) on $\mathbb{R}^{n+1}$ such that 
 $$\lim_{t \rightarrow 0^{-}} \|u\|_{L^{2+\delta}_x(B_1 \times \{-t\})} = \infty, \quad \lim_{t \rightarrow 0^{-}} \|\nabla u\|_{L^{2+\delta}(B_1 \times (-1,\,-t))} = \infty.$$
\end{cor}
\noindent (The ellipticity ratio $\lambda / \Lambda$ degenerates as $\delta \rightarrow 0$, in accordance with the higher-integrability results).
We conclude, as in the elliptic case, that solutions to parabolic systems are only slightly better than their energy class.

Our remaining results are Liouville theorems for (\ref{SelfSimEquation}).
It is natural to ask whether one can construct solutions that decay any faster than than we managed. Our first Liouville theorem shows this is not possible:
\begin{thm}\label{Liouville}
 Assume that $w \in H^1_{loc}(\mathbb{R}^n)$ solves (\ref{SelfSimEquation}), with $|w|  = O(|x|^{-\mu})$ and $2\mu \geq n$. Then $w \equiv 0$.
\end{thm}
\noindent There are nontrivial $-\mu$-homogeneous solutions to elliptic systems of the form $\text{div}(A(x)\nabla u) = 0$ in $\mathbb{R}^n$ provided $2\mu < n-2$, and there is a Liouville theorem for $-\mu$-homogeneous solutions 
on $\mathbb{R}^n \backslash \{0\}$ in the equality case (see \cite{M2}).  Thus, Theorems \ref{Counterexample} and \ref{Liouville} mirror the elliptic results in dimension $n+2$.
This agrees with the observation that the parabolic energy $L^{\infty}_t(L^2_x) + L^2_t(H^1_x)$ in $\mathbb{R}^{n+1}$ and the elliptic energy $H^1$ in $\mathbb{R}^{n+2}$ are invariant under the matching rescalings
$$u \rightarrow \lambda^{n/2}u(\lambda x,\,\lambda^2 t), \quad \text{ resp. } \quad u \rightarrow \lambda^{n/2}u(\lambda x).$$

Theorem \ref{Liouville} is a consequence of parabolic energy estimates. We can extend it to the ``elliptic regime" $2\mu \geq n-2$ when $w$ has the monotonicity property
\begin{equation}\label{Monotonicity}
 (2\mu + x \cdot \nabla)|w|^2 \geq 0:
\end{equation}
\begin{thm}\label{EllipticLiouville}
 Assume that $w \in H^1_{loc}(\mathbb{R}^n)$ solves (\ref{SelfSimEquation}), with $|w| = O(|x|^{-\mu})$ and $2\mu \geq n-2$. If in addition
 $w$ satisfies (\ref{Monotonicity}), then $w \equiv 0$.
\end{thm}
\noindent It is easy to check that previous examples (\cite{FM}, and \cite{SJM}, \cite{SJ2} for more general systems) satisfy condition (\ref{Monotonicity}), which explains
why they have ``elliptic" discontinuities (that is, $n \geq 3$ and $2\mu < n-2$).

\vspace{2mm}

The paper is organized as follows. In Section \ref{CounterexampleSection} we prove Theorem \ref{Counterexample}. In Section \ref{LiouvilleTheorems} we prove
Theorems \ref{Liouville} and \ref{EllipticLiouville}. Finally, in Section \ref{OpenQuestions} we list some open questions.

\section{Proof of Theorem \ref{Counterexample}}\label{CounterexampleSection}

In this section we prove Theorem \ref{Counterexample}. We exploit the useful observation that if $\text{Im}(A)$ is symmetric, then the ellipticity
condition (\ref{Ellipticity}) is satisfied provided $\text{Re}(A)$ is uniformly positive definite and $|A|$ is bounded (see \cite{F}). 

\begin{rem}
 Heuristically, this structure allows strong coupling between components when we view (\ref{Equation}) as the system (\ref{System}).
 The example in \cite{M1} has skew-symmetric imaginary coefficients, which corresponds to the symmetry $B^{kl}_{\alpha\beta} = B^{lk}_{\beta\alpha}$ of 
 the system coefficients. In that case it is important to estimate the size of $\text{Im}(A)$ since it affects the ellipticity condition.
\end{rem}


\subsection{Reduction to ODE System}
We first reduce (\ref{SelfSimEquation}) to an ODE system. Let $r = |x|$ and let $\nu = r^{-1}x$ be the unit radial vector.
We search for solutions of the form
\begin{equation}\label{SolutionForm}
w = \varphi(r)g(\nu)e^{-i\log r}.
\end{equation}
Then
\begin{equation}\label{Gradient}
\nabla w = ge^{-i\log r}(\varphi'(r) - ir^{-1}\varphi)\nu + \varphi(r)e^{-i\log r}r^{-1}\nabla_{S^{n-1}}g.
\end{equation}
Here and below $\nabla_{S^{n-1}}$ and $\Delta_{S^{n-1}}$ denote the usual gradient and Laplace operators on the sphere. If 
$$B = f(r) \nu \otimes \nu + h(r)(I - \nu \otimes \nu)$$ 
then we have
$$B\nabla w = ge^{-i\log r}r^{n-1}(f\varphi' - ir^{-1}f\varphi) \frac{\nu}{r^{n-1}} + h\varphi e^{-i\log r}r^{-1} \nabla_{S^{n-1}}g.$$
We will choose $\varphi$ such that $\varphi'$ and $r^{-1}\varphi$ are bounded. Using that $\nu / r^{n-1}$ is divergence-free away from the origin we compute
\begin{align*}
\text{div}(B\nabla w) &=  \\
&ge^{-i\log r}\left[\frac{(r^{n-1}f\varphi')'}{r^{n-1}} - \left(f - \frac{\Delta_{S^{n-1}}g}{g} h\right)\frac{\varphi}{r^2} 
- i\left(\frac{(r^{n-2}f \varphi)'}{r^{n-1}} + \frac{f \varphi'}{r}\right)\right].
\end{align*}
Let $g$ be an eigenfunction of $\Delta_{S^{n-1}}$ with eigenvalue $-\lambda_g < 0$. Then the previous expression becomes
$$\text{div}(B\nabla w) = ge^{-i\log r} \left[\frac{(r^{n-1}f\varphi')'}{r^{n-1}} - (f + \lambda_g h)\frac{\varphi}{r^2} 
- i \frac{(r^{n-2}f \varphi^2)'}{r^{n-1}\varphi}\right].$$
Thus, if we take coefficients
\begin{equation}\label{Coefficients}
 A = \alpha I + i (\beta(r) \nu \otimes \nu + \gamma(r)(I - \nu \otimes \nu))
\end{equation}
with $\alpha > 0$ constant, and $g$ is any linear function restricted to the sphere, we obtain
\begin{align*}
 \text{div}(A\nabla w) &=  ge^{-i\log r} \left[\alpha \left(\frac{(r^{n-1}\varphi')'}{r^{n-1}} - n\frac{\varphi}{r^2}\right) + \frac{(r^{n-2}\beta\varphi^2)'}{r^{n-1} \varphi} \right. \\
&+ \left. i \left(\frac{(r^{n-1}\beta\varphi')'}{r^{n-1}} - (\beta + (n-1)\gamma) \frac{\varphi}{r^2} - \alpha \frac{(r^{n-2}\varphi^2)'}{r^{n-1}\varphi}\right)\right].
\end{align*}
Since
$$iw + \mu w + x \cdot \nabla w = ge^{-i\log r}(\mu \varphi + r\varphi'),$$
the equation (\ref{SelfSimEquation}) becomes the ODE system
\begin{equation}\label{ODESystem}
\begin{cases}
\frac{(r^{n-2}\beta\varphi^2)'}{r^{n-1} \varphi} = \frac{1}{2}(\mu \varphi + r\varphi') + n\alpha \frac{\varphi}{r^2} - \alpha \frac{(r^{n-1}\varphi')'}{r^{n-1}}, \\
(n-1)\gamma \frac{\varphi}{r^2} = -\alpha \frac{(r^{n-2}\varphi^2)'}{r^{n-1}\varphi} + \frac{(r^{n-1}\beta\varphi')'}{r^{n-1}} - \beta \frac{\varphi}{r^2}.
\end{cases}
\end{equation}

We will fix $\varphi \sim r^{-\mu}$ and $\alpha > 0$ depending on $\mu$. Then the first equation determines $\beta$, and the second one $\gamma$.
By the remark at the beginning of the section, the point is to make choices such that $\beta$ and $\gamma$ are bounded.

\subsection{Solving the ODE System}
Integrating the first equation in (\ref{ODESystem}) we obtain
\begin{equation}\label{IntegratedEquation}
\begin{split}
\beta &= \frac{1}{4}\left(r^2 + \frac{2\mu-n}{r^{n-2}\varphi^2} \int_0^r \varphi^2(s)s^{n-1}\,ds \right) \\
&+ \frac{n\alpha}{r^{n-2}\varphi^2} \int_0^r \varphi^2(s)s^{n-3}\,ds \\
&+ \frac{\alpha}{r^{n-2}\varphi^2} \int_0^r \varphi'^2(s)s^{n-1}\,ds - \alpha\frac{r\varphi'}{\varphi}.
\end{split}
\end{equation}

\begin{rem}
It follows easily that if $2\mu \geq n$ and $\varphi = O(r^{-\mu})$, then $\beta$ is unbounded (compare to Theorem \ref{Liouville}).
\end{rem}

We define
\begin{equation}
 \varphi(r) = \begin{cases}
            r,\, \quad 0 \leq r < 3/4 \\
            r^{-\mu} + C_{\mu}r^{-\mu-2}, \quad r > 1 \\
            \text{positive and smooth,} \quad 1/2 < r < 3/2
           \end{cases}
\end{equation}
where $C_{\mu} \geq 0$ will be chosen later. 

\begin{rem}\label{InterestingCases}
 By Theorem \ref{EllipticLiouville} it will be necessary to take $C_{\mu} > 0$ when $2\mu \geq n-2$ (and in particular, to generate discontinuities in the case $n = 2$).
\end{rem}

\vspace{2mm}

For $r < 3/4$ it is easy to check that $\beta$ and $\gamma$ are of the form $c_1(n,\,\alpha) + c_2(n,\,\mu)r^2$ (with $c_i$
linear in $\alpha$ and $\mu$) so we only need to analyze the solutions for $r$ large. We divide into three cases.

\vspace{2mm}

{\bf Case $1$: $2\mu < n-2$.} We take $C_{\mu} = 0$ and $\alpha = 1$. It is easy to check that $\beta$ and $\gamma$ have 
the form $c_1 + c_2r^{2-n+2\mu}$ for $r > 1$, which is bounded.

\vspace{2mm}

{\bf Case $2$: $n-2 < 2\mu < n$.} Now the quantities
$$D := \int_0^{\infty} (\varphi^2 - s^{-2\mu})s^{n-1}\,ds, \quad E := \int_0^{\infty} \varphi^2 s^{n-3}\,ds, \quad F := \int_{0}^{\infty} \varphi'^2 s^{n-1}\,ds$$
are bounded, for any fixed $C_{\mu} \geq 0$. The solution (\ref{IntegratedEquation}) becomes
$$\beta = \left(-\frac{n-2\mu}{4}D + \alpha(nE + F)\right)r^{2\mu - n + 2} + \mathcal{R}(1).$$
Here and below, $\mathcal{R}(1)$ denotes any smooth function on $(1,\,\infty)$ whose $j^{th}$ derivative is $O(r^{-j})$ as $r \rightarrow \infty$ 
for each $j \geq 0$. Using the definition of $\varphi$ we estimate
\begin{align*}
 D &\geq -\int_0^1 s^{n-1-2\mu}\,ds + 2C_{\mu}\int_1^{\infty} s^{n-3-2\mu}\,ds \\
 &= -\frac{1}{n-2\mu} + \frac{2C_{\mu}}{2\mu - n + 2}.
\end{align*}
We conclude that
$$-\frac{n-2\mu}{4}D \leq \frac{1}{4} - \frac{n-2\mu}{2(2\mu-n+2)}C_{\mu} < 0$$
provided we choose $C_{\mu}$ large. We may then choose $\alpha > 0$ small so that
$$-\frac{n-2\mu}{4}D + \alpha(nE + F) = 0,$$
hence
$$\beta = \mathcal{R}(1).$$
Solving the second equation in (\ref{ODESystem}) for $\gamma$ gives
$$\gamma = \mathcal{R}(1),$$
which completes this case.

\vspace{2mm}

{\bf Case $3$: $2\mu = n-2$.} This case is similar to the case $2\mu > n-2$, except to leading order $\beta$ grows logarithmically. Computing
(\ref{IntegratedEquation}) gives
$$\beta = \left(-C_{\mu} + \alpha \left(n + \frac{1}{4}(n-2)^2\right)\right)\log r + \mathcal{R}(1).$$
Choosing $C_{\mu}$ and $\alpha$ to satisfy the relation
$$C_{\mu} = \left(n + \frac{1}{4}(n-2)^2\right)\alpha$$
we arrive at the same conclusion as in Case $2$, completing the construction.

\subsection{Proof of Theorem \ref{Counterexample}}
\begin{proof}[{\bf Proof of Theorem \ref{Counterexample}}:]
 For $0 \leq 2\mu < n$, take $\varphi,\, g,\, \alpha,\, \beta,\, \gamma$ as constructed above. Then the function
 $$w = \varphi(r)g(\nu)e^{-i \log r}$$
 solves the equation (\ref{SelfSimEquation}) in $\mathbb{R}^n$ with bounded coefficients
 $$A = \alpha I + i(\beta(r) \nu \otimes \nu + \gamma(r) (I-\nu \otimes \nu))$$
 and has the asymptotics (\ref{Asymptotics}). Since $\alpha > 0$ is constant and $\text{Im}(A)$ is symmetric, 
 the coefficients satisfy the ellipticity condition (\ref{Ellipticity}), completing the proof.
\end{proof}

\begin{rem}
In our construction, $w$ is Lipschitz but no better at $0$, and smooth but not analytic away from $0$. This is a consequence
of choices we made for computational convenience. It is not hard to modify the construction so that $w$ is analytic on $\mathbb{R}^n$, e.g.
by taking $w = \varphi(r)g(\nu)e^{-\frac{i}{2}\log(1+r^2)}$ with $g$ as above and
$$\varphi = r\left((1+r^2)^{-\frac{\mu+1}{2}} + C_{\mu}(1+r^2)^{-\frac{\mu+3}{2}}\right).$$ 
The coefficients $A(x)$ also become analytic with these modifications. 

\end{rem}

\section{Liouville Theorems}\label{LiouvilleTheorems}

In this section we prove the Liouville theorems Theorem \ref{Liouville} and Theorem \ref{EllipticLiouville}.

\subsection{Proof of Theorem \ref{Liouville}}
\begin{proof}[{\bf Proof of Theorem \ref{Liouville}}]
Let $\psi \in C^{\infty}_0(\mathbb{R}^n)$ be real-valued. Multiplying (\ref{SelfSimEquation}) by $\overline{w}\psi^2$ we obtain
\begin{equation}\label{KeyInequality}
2\text{Re}\left( \text{div}(A\nabla w)\overline{w}\psi^2 \right) = \frac{1}{2}(2\mu |w|^2 + x  \cdot \nabla |w|^2)\psi^2.
\end{equation}
Integrating by parts and using the ellipticity condition (\ref{Ellipticity}) we get
\begin{equation}\label{Caccioppoli}
\begin{split}
\int_{\mathbb{R}^n} (-\lambda |\nabla w|^2\psi^2 &+ C(\lambda,\,\Lambda)|w|^2|\nabla \psi|^2)\,dx \\
&\geq \frac{2\mu - n}{2} \int_{\mathbb{R}^n} |w|^2\psi^2\,dx - \frac{1}{2}\int_{\mathbb{R}^n} |w|^2 x \cdot \nabla (\psi^2)\,dx.
\end{split}
\end{equation}
Since $2\mu \geq n$, the first term on the right side is non-negative.

We now fix our choice of $\psi$. Let $\psi_1$ be a smooth, radially decreasing function supported in $B_2$ with $\psi_1 \equiv 1$ in $B_1$, 
and let $\psi_R := \psi_1(R^{-1}x)$. Take $\psi = \psi_R$. Then the second term on the right side of (\ref{Caccioppoli}) is non-negative, so the
right side is non-negative. Using that $|w|^2|\nabla \psi|^2 = O(R^{-2\mu-2})$ in $B_{2R} \backslash B_R$ we conclude that
$$\int_{B_R} |\nabla w|^2\,dx = O(R^{n-2\mu-2}) = O(R^{-2}),$$
completing the proof.
\end{proof}

\subsection{Proof of Theorem \ref{EllipticLiouville}}

\begin{proof}[{\bf Proof of Theorem \ref{EllipticLiouville}}]
We start again with the identity (\ref{KeyInequality}).
By (\ref{Monotonicity}) the right side of (\ref{KeyInequality}) is non-negative. Integrating by parts gives the Caccioppoli inequality
$$\int_{\mathbb{R}^n} |\nabla w|^2\psi^2\,dx \leq C(\lambda,\,\Lambda) \int_{\mathbb{R}^n} |w|^2|\nabla \psi|^2\,dx.$$
Choosing $\psi$ as before, we recover the inequality
$$\int_{B_{R}} |\nabla w|^2\,dx = O(R^{n-2\mu-2}),$$
which proves the theorem when $2\mu > n-2$.
In the critical case $2\mu = n-2$, use instead 
$$\psi = \begin{cases}
1 \text{ in } B_1, \\
1-\log(r)/\log(R) \text{ in } B_R \backslash B_1, \\
0 \text{ in } \mathbb{R}^n \backslash B_R
\end{cases} $$
to obtain
$$\int_{B_{\sqrt{R}}} |\nabla w|^2\,dx = O\left(\frac{1}{\log R}\right).$$
\end{proof}

\section{Some Questions}\label{OpenQuestions}

To conclude we list some open questions.

\vspace{2mm}

\begin{enumerate}
 \item Our examples have coefficients with symmetric imaginary part. 
 Similar constructions might be possible with skew-symmetric imaginary coefficients, using techniques from \cite{M1}.
 In this setting the imaginary coefficients play a role in ellipticity. 
 \vspace{2mm}

 \item For elliptic systems there is a sharp condition on the spectrum of the coefficients that guarantees continuity of solutions \cite{Ko}. 
 Sufficient conditions are known in the parabolic case (\cite{Ko}, \cite{Ka}). It would be interesting to investigate how closely our
 counterexamples match these conditions.

 \vspace{2mm}

 \item Solutions to parabolic systems in dimension $n \geq 3$ can be discontinuous on very large sets \cite{SJ1}. 
 It is natural to ask how large the discontinuity set can be when $n = 2$. 
 Known results imply spatial continuity at almost every time, which is false when $n \geq 3$ by elliptic examples.

 \vspace{2mm}

 \item Parabolic systems with the quasilinear structure 
 \begin{equation}\label{QuasilinearSystem}
 u_t - \text{div}(A(u)\nabla u) = 0
 \end{equation}
 have a well-developed partial regularity theory and are important in applications \cite{GS}. Here the coefficients depend smoothly on $u$.
Constructing solutions to (\ref{QuasilinearSystem}) becomes easier when $u \in \mathbb{R}^m$ for $m$ large because there is more room to ``disperse $u$.'' 
 \cite{M1} contains examples of discontinuity formation for (\ref{QuasilinearSystem}) when $n = 2,\, m = 4$.
 One can improve to $n = 2,\, m = 3$ using similar techniques \cite{M3}. 
 Continuity for solutions to (\ref{QuasilinearSystem}) in the case $n = m = 2$ (in particular, the $\mathbb{C}$-valued scalar case) remains open. 
 It seems possible in view of Theorem \ref{EllipticLiouville} that the restrictive geometry of the target could play in favor of regularity (see the discussion in \cite{M3}).

\end{enumerate}



\section*{acknowledgements}
This work was supported by NSF grant DMS-1501152 and by the ERC grant ``Regularity and Stability in PDEs."
I am grateful to John Ball and Jan Kristensen for discussions, and for the generous hospitality of the Oxford Mathematical Institute
during the time this work was completed.






\end{document}